\documentclass[a4paper,11pt]{article}

\usepackage{MyLaTeX}
\usepackage{multirow}

\numberwithin{defn}{section}
\input xy
\xyoption{all}
\newcommand{\Supp}{\operatorname{Supp}}

\title{On the Unit Conjecture for Supersoluble Group Rings, I}
\author{David A.~Craven and Peter Pappas}
\date{September 2010}

\begin{document}
\maketitle

\begin{abstract}We introduce structure theorems for the study of the unit conjecture for supersoluble group rings and apply our results to the (Passman) fours group \[ \Gamma=\Gen{x,y}{y^{-1}x^2y=x^{-2},\;x^{-1}y^2x=y^{-2}}.\] We show that over any field $K$, the group algebra $K\Gamma$ has no non-trivial units $\sigma$ of length $L(\sigma )\leq 3$, and find that the Promislow set can never be the support of a unit in $K\Gamma$. We conclude our work with an introduction to the theory of {\it consistent chains} toward a preliminary analysis of units of higher length in $K\Gamma$.\end{abstract}

\section{Introduction}

The unit conjecture for group algebras asserts that if $K$ is a field and if $G$ is a torsion-free group, then every unit\footnote{\emph{Unit} here means two-sided unit. If $K$ is a field of characteristic $0$, then Kaplansky's theorem \cite[p.\ 38]{miliessehgal} states that every unit in $KG$ is two-sided. If $K$ has characteristic $p$ and $G$ is polycyclic-by-finite, then Farkas--Marciniak obtain a similar result using a Witt ring construction \cite{farkasmarc1982}. The general result for group algebras over fields of characteristic $p$ remains open.} of the group algebra $KG$ is \emph{trivial}; that is, every unit is of the form $\lambda g$ for some $\lambda\in K\setminus\{0\}$ and $g\in G$ \cite{miliessehgal} \cite{passman} \cite{sehgal} \cite{sehgal2}. The best result to date is entirely group-theoretic, concerning group algebras of unique-product groups \cite{passman} \cite{passman2}. (A group $G$ is said to be a \emph{unique-product group} if, given any two non-empty finite subsets $X$ and $Y$ of $G$, there exists an element $g\in G$ having a unique representation of the form $g=xy$  with $x\in X$ and $y\in Y$.) Unique-product groups typify ordered, right-ordered, locally indicable groups and for some time it remained an open question whether there exist torsion-free groups that are not unique-product groups. Using small cancellation theory, Rips and Segev \cite{ripssegev1987} gave the first example of a torsion-free group that is not a unique-product group.

For the unit conjecture beyond unique-product groups, it is clear that one should consider finitely generated, torsion-free, abelian-by-finite groups; that is, groups with a short exact sequence
\[ 1\to A\to G\to G/A\to 1\]
with $A$ abelian and $G/A$ finite. If $G/A$ is cyclic then $G$ is right-orderable, and therefore a unique-product group, so nothing new occurs. The simplest example where $G/A$ is non-cyclic is
\[ \Gamma=\gen{x,y\mid x^{-1}y^2x=y^{-2},\;y^{-1}x^2y=x^{-2}},\]
which satisfies the short exact sequence
\[ 1\to \Z^3\to \Gamma\to \Z/2\Z\times \Z/2\Z\to 1.\]
Called the `fours group', $\Gamma$ was introduced by Passman \cite[p.\ 606]{passman} and shown to be torsion-free and non-right orderable. Promislow \cite{promislow1988}, using a random search algorithm, exhibited a 14-element subset $\ms P\subs \Gamma$ such that $\ms P\cdot \ms P$ has no unique product\footnote{It is an open question as to whether every unique-product group is right orderable.}. Since then, very little progress on the unit conjecture has been made, and it has been an open question whether the Promislow set $\ms P$ could be the support of a unit over some field $K$.

In this paper we show that the answer is `no'. To obtain our result we first derive a \emph{splitting theorem} for units in $K\Gamma$. This is implicit in earlier work of Cohn \cite{cohn1968} and Lewin \cite{lewin1972}, and is a direct consequence of Passman's work \cite[Theorem 13.3.7]{passman}. The group $\Gamma$ is supersoluble and contains a normal subgroup $N$ such that $\Gamma/N$ is infinite dihedral. This leads to a length function $L:K\Gamma\to \N\cup\{-\infty\}$ and we show, via the splitting theorem, that if $u\in K\Gamma$ is a unit then $L(u)=L(u^{-1})$. On the other hand, the group $\Gamma$ being abelian-by-finite, with $A=\Z^3$ in the notation above, induces a faithful representation $\eta:K\Gamma\hookrightarrow M_4(KA)$, and we find, for $\alpha\in K\Gamma$, that $\alpha$ is a unit of $K\Gamma$ if and only if $\det(\eta(\alpha))$ is a non-zero element of the field\footnote{This result is known more generally for crystallographic groups (though to the best of our knowledge unpublished). We thank Dan Farkas for conveying this fact to us; out proof is elementary and we include it for completeness.}. Our main result then shows that there are no non-trivial units in $K\Gamma$ of length at most $3$. Applying a specific automorphism of $K\Gamma$ allows us to show that the Promislow set $\ms P$ can never be the support of a unit in $K\Gamma$ for any field $K$.

We conclude with a discussion of how our techniques apply to the higher-length situation, which is the subject of the sequel to this paper. To this end we introduce the theory of {\it consistent chains} toward a preliminary analysis of units of higher length in $K\Gamma$.

\section{A Splitting Theorem for Supersoluble Groups}
\label{splitsec}

Let $G$ be a group, and assume that $N$ is a normal subgroup of $G$ such that $G/N$ is the infinite dihedral group, generated by involutions $Nx$ and $Ny$. Write $X=\gen{N,x}$ and $Y=\gen{N,y}$. Let $W$ be the set of all alternating words in $x$ and $y$. For example, $xyxy$ is an element of $W$, and we say that it \emph{starts} in $X$ and \emph{ends} in $Y$. Since $G/N=(X/N)\ast (Y/N)$, it follows from \cite[Theorem 9.2.9]{passman} that $W$ is a transversal for $N$ in $G$. If $g\in Nw$, then we let the starting and ending properties of $w$ carry over to $g$.

We now define a length function on $KG$. The \emph{length} of a word $w\in W$, denoted by $L(w)$, is the number of factors that occur in it; the empty word, $w=1$, has length $0$, and the example $xyxy$ has length $4$. We extend the length function $L$ in two ways: firstly, if $g\in G$ then there exists a unique $w\in W$ with $g\in Nw$, and we define $L(g)=L(w)$; and secondly, if $\alpha\in KG$ with $\alpha$ non-zero, then we set $L(\alpha)$ to be the maximum of $L(g)$, where $g\in \Supp \alpha$. Finally, set $L(0)$ to be $-\infty$. From $W\subs G\subs KG$, we see that the definition of $L$ is consistent. 

For $w_1,w_2\in W$  we say that the product $w_1w_2$ is \emph{non-overlapping} if no cancellation occurs.  In this case, 
\[ L(w_1w_2)=L(w_1)+L(w_2).\]
On the other hand, if the product $w_1w_2$ overlaps, then $L(w_1w_2)$ is strictly less than $L(w_1)+L(w_2)$. In this case, if $w_1$ ends in $X$ (and hence $w_2$ starts in $X$) then we say that the overlap is in $X$, and similarly for overlapping in $Y$. 

With these assumptions and notation, we can now state our first result, which is a direct consequence of the work of Cohn \cite{cohn1968} and Lewin \cite{lewin1972}. The proof that we give follows that given in \cite[Theorem 13.3.7]{passman}.

\begin{thm}\label{splittingthm}Let $K$ be a field and let $G$ be a group with a normal subgroup $N$ as above. Assume that $KG$ has no proper divisors of zero and that $KN$ is an Ore domain. Suppose that for some $\sigma$, $\tau\in KG\setminus\{0\}$ we have that $\sigma\tau\in KN$. There exist $\alpha_1,\dots,\alpha_s$, $\beta_1,\dots,\beta_s\in KN$, $\gamma_1,\dots,\gamma_s\in\{x,y\}$ with $L(\gamma_1\ldots\gamma_s)=s$, such that
\[ \sigma(\alpha_1+\beta_1\gamma_1)\ldots(\alpha_s+\beta_s\gamma_s)\in KN\setminus\{0\}.\]
\end{thm}
\begin{pf}Assume that $\sigma$ and $\tau$ are non-zero elements of $KG$ with $\sigma\tau\in KN$. Then $L(\sigma)$, $L(\tau)\geq 0$ and, moreover, $\sigma\tau$ is non-zero. We prove the theorem by induction on $L(\sigma)$.

If $L(\sigma)=0$ then $\sigma\in KN$, so that $\sigma(1+0\cdot x)\in KN\setminus\{0\}$ yields the desired result; therefore, we may assume that $L(\sigma)>0$, and by induction the result holds for all such $\bar\sigma$ and $\bar\tau$ with $L(\bar\sigma)<L(\sigma)$. Since $L(\sigma)>0$, we see that $\sigma$ is not in $KN$, and so therefore neither is $\tau$: hence $L(\tau)>0$. Let $L(\sigma)=m$ and $L(\tau)=n$. We proceed in a series of three steps, the first two of which are exactly those given in the proof of \cite[Theorem 13.3.7]{passman}. Because of this, we will suppress the proofs of the first two steps, and invite the interested reader to consult \cite{passman}.

\medskip

\noindent {\bf Step 1:} \emph{The products of maximal-length elements overlap in the same group.}

\medskip

We assume, by symmetry, that the products of maximal-length elements overlap in $X$. Write $\sigma=\sigma'+\sigma''$, where $\Supp \sigma'$ is given by all those elements $g\in\Supp\sigma$ with either $L(g)=L(\sigma)=m$ or with $L(g)=m-1$ and with $g$ ending in $Y$. All elements of length $m$ in $\Supp\sigma$ end in $X$ so that $\sigma'=\sum a_i\ep_i$, where the elements $a_i$ of $W$ have length $m-1$ and end in $Y$, and where $\ep_i\in KX\setminus \{0\}$. Similarly, write $\tau=\tau'+\tau''$,  where $\Supp \tau'$ consists of all those elements $g\in\Supp\tau$ with either $L(g)=L(\tau)=n$, or with $L(g)=n-1$ and with $g$ starting in $Y$. It follows that $\tau'=\sum \delta_jb_j$, where the elements $b_j$ come from $W$ all have length $n-1$ and start in $Y$, and where $\delta_j\in KX\setminus\{0\}$.

\medskip

\noindent {\bf Step 2:} \emph{The products $\ep_i\delta_j$ all belong to $KN$.} See Step 2 of \cite[Theorem 13.3.7]{passman}.

\medskip

\noindent {\bf Step 3:} \emph{The inductive step.}

Since $N\normal G$, \cite[Lemma 13.3.5(ii)]{passman} implies that the set $T=KN\setminus\{0\}$ of regular elements of $KN$ is a right divisor set of regular elements of $KG$. Now $\ep_1\delta_1\in T$ and $\ep,\tau\in KG$, so there exist elements $\eta\in T$ and $\rho\in KG$ with
\[ (\ep_1\delta_1)\rho=(\ep_1\tau)\eta.\]
Thus, because $\ep_1$ and $\eta$ are regular elements of $KG$ and $\tau$ is non-zero, we conclude that $\rho\neq0$ and $\delta_1\rho=\tau\eta$. This yields
\[ (\sigma\delta_1)\rho=(\sigma\tau)\eta\in KN,\]
so that $(\sigma\delta_1)\rho\in KN$.

We now compute the length of $(\sigma\delta_1)\rho$. We observe that $\sigma\delta_1\neq 0$ since $\sigma\neq 0$, and $\delta_1\neq 0$ implies that $\delta_1$ is not a zero divisor in $KG$. Thus $L(\sigma\delta_1)\geq 0$. Moreover
\[ \sigma'\delta_1=\sum_i x_i(\ep_i\delta_1),\]
and since $L(x_i)=m-1$ and $\ep_i\delta_1\in KN$, by Step 2, we conclude that $L(\sigma'\delta_1)\leq m-1$. Since $L(\sigma'')\leq m-1$ and $\delta_1\in KX$, we have
\[ L(\sigma''\delta_1)\leq L(\sigma'')+L(\delta_1)\leq (m-1)+1=m.\]
If equality occurs then there exist elements $g\in \Supp \sigma''$, $h\in\Supp \sigma'$ with $L(g)=m-1$, $L(h)=1$, and with $gh$ non-overlapping. However, $L(g)=m-1$, and $g\in \Supp \sigma''$ implies that $g$ ends in $X$ and $h$ starts in $X$. Therefore, the product does overlap, and this case cannot occur. Hence $L(\sigma''\delta_1)\leq m-1$, and from $\sigma\delta_1=\sigma'\delta+\sigma''\delta_1$, it follows that
\[ 0\leq L(\sigma\delta_1)\leq m-1<L(\sigma).\]

By induction, there exist $\alpha_1,\dots,\alpha_s$, $\beta_1,\dots,\beta_s\in KN$, $\gamma_1,\dots,\gamma_s\in\{x,y\}$ with $L(\gamma_1\ldots\gamma_s)=s$, such that
\[ \sigma(\alpha_1+\beta_1\gamma_1)\ldots(\alpha_s+\beta_s\gamma_s)\in KN\setminus\{0\}.\]
The result now follows, noting that $\delta_1=\alpha+\beta a\neq 0$ for some $\alpha,\beta\in KN$.\end{pf}

This means that if $\sigma\tau=1$ then we may write $\tau$ as a product $t$ of linear terms (i.e., $\alpha_i+\beta_i\gamma_i$ with $\alpha_i,\beta_i\in KN$ and $\gamma_i\in\{x,y\}$) times the \emph{inverse} of some element $\ep\in KN$. Either we get $\sigma t=\ep$ or, by formally inverting the elements of $KN$, $\sigma t\ep^{-1}=1$. We will refer to this product as a \emph{splitting} for $\tau$. Note that this splitting is not unique in general; we will discuss this problem later. We will tend to write $\sigma=\eta^{-1}s$ for a splitting of $\sigma$ and $\tau=t\ep^{-1}$ for a splitting of $\tau$. Of course, since all units of $KG$ are two-sided, $\sigma\tau=1$ implies $\tau\sigma=1$, so we may get a splitting $\sigma=s\eta^{-1}$ for some (potentially different) $s$ and $\eta$, and similarly for $\tau$.

\section{Using the Splitting Theorem}
\label{sec:conseqsplitting}

The splitting theorem of the previous section is a powerful tool for analyzing units in supersoluble groups. If we analyze a `minimal' counterexample $G$ to the unit conjecture, we may assume that all subgroups of $G$ of smaller Hirsch length satisfy the unit conjecture over a given field $K$; we call such a group a \emph{UC-proper} group. Our first theorem gives information on the inverse of a unit, and the second gives information on the structure of words of maximal length in $\sigma$.

\begin{thm}\label{thm:samelength} If $\sigma,\tau\in KG\setminus \{0\}$ such that $\sigma\tau\in KN$, then $L(\sigma)=L(\tau)$.\end{thm}
\begin{pf}In the notation of Step 3, we have $\delta_1\rho=\tau\eta$, and by Theorem \ref{splittingthm},
\[ \delta_1\rho=(\alpha_1+\beta_1\gamma_1)\ldots(\alpha_s+\beta_s\gamma_s)\]
with $L(\delta_1\rho)=s$. But $L(\delta_1\rho)=L(\tau\eta)$, so that $s=L(\delta_1\rho)=L(\tau)$. Observe that the argument in Theorem \ref{splittingthm} is left-right symmetric. Let $\tau'=\delta_1\rho$; we have $\sigma\tau'\in KN\setminus\{0\}$.

Proceeding as in Steps 1 and 2, and using $T=KN\setminus\{0\}$ as a left divisor set of regular elements of $KG$, we get $\ep_1\delta_1\in T$ and $\sigma\delta_1\in KG$, so that there exist elements $\eta'\in T$ and $\rho'\in KG$ with $\rho'\ep_1\delta_1=\eta'\sigma\delta_1$, and as before we conclude that $\rho'\ep_1=\eta'\sigma$. Thus $L(\rho'\ep_1)=L(\sigma)=t$. An inductive argument yields
\[ \sigma'=\rho'\ep_1=(\alpha'_1+\beta'_1\gamma'_1)\ldots(\alpha'_t+\beta'_t\gamma'_t)\]
with $\alpha_i',\beta_j'\in KN$ and $L(\gamma_1'\ldots\gamma_t')=t$. Hence,
\[\left[(\alpha'_1+\beta'_1\gamma'_1)\ldots(\alpha'_t+\beta'_t\gamma'_t)\right]\left[\vphantom{\alpha'_1}(\alpha_1+\beta_1\gamma_1)\ldots(\alpha_s+\beta_s\gamma_s)\right]\in KN\setminus\{0\}.\]

Observe that $\gamma_1'\ldots\gamma_t'$ and $\gamma_1\ldots\gamma_s$ are the unique words in $\sigma'$ and $\tau'$ of maximal length. By our remarks in Theorem \ref{splittingthm}, the elements $\gamma_t'$ and $\gamma_1$ belong to the same group, say $X$. If $(\alpha_t'+\beta_t'a)(\alpha_1+\beta_1 a)$  does not lie in $KN$, then this contains some term of the form $\nu x$. Arguing as in Step 2 shows that
\[ \gamma_1'\ldots \gamma_{t-1}'x\gamma_2\ldots \gamma_s\]
would occur only once in the product $\sigma'\tau'$, which is impossible, since this must be cancelled off. Thus $(\alpha_t'+\beta_t'x)(\alpha_1+\beta_1 x)\in KN$, so by induction $t-1=s-1$. Thus $s=t$ as desired.
\end{pf}

\begin{cor}\label{oneeltmaxlength} Suppose that $\sigma\tau=1$. Then there is only one word of maximal length in $\sigma$. If $\sigma=\sigma^*$, then $L(\sigma)$ is odd; i.e., the word of maximal length in $\sigma$ starts and ends in the same group.\end{cor}
\begin{pf} By Step 1, the products of maximal-length words in $\sigma$ and $\tau$ all overlap in the same group; thus $\sigma$ has only one maximal-length word. If $\sigma=\sigma^*$, this must begin and end in the same group, and so has odd length.
\end{pf}

We now want to analyze the element $\eta$ of $KN$ that we invert to go from the split form of $\sigma$ to $\sigma$ itself. As in the previous section, write $W$ for the set of all words in $x$ and $y$, creating a transversal to $N$ in $G$. For a given element $\sigma\in KG$, let $I$ denote the subset of all words in $W$ in the support of $\sigma$.

\begin{prop}\label{prop:leftgcdis1} Let $G$ be a UC-proper, supersoluble group and let $\sigma$ be a non-trivial unit. Write
\[ \sigma=\sum_{w \in I} a_ww,\]
where $a_w\in KN$. The left-gcd of the $a_w$ is 1. In other words, if $\sigma=\ep \sigma'$ with $\ep\in KN$ then $\ep=\lambda g$ for $\lambda\in K$ and $g\in N$.
\end{prop}
\begin{pf} If $\sigma=\ep\sigma'$ is a unit, then $\sigma\tau=\ep\sigma'\tau=1$, so that $\ep$ is a unit. Since $G$ is UC-proper, $\ep$ is a trivial unit, as claimed.\end{pf}

If $\sigma$ is a unit and we write $\sigma=\eta^{-1}s$, where $s$ is a split, by the previous proposition we must have that the $\eta^{-1}$ must cancel off the entire gcd of the coefficients in front of the words in $I$.

\begin{cor} Let $G$ be a UC-proper, supersoluble group, and let $\sigma$ be a non-trivial unit, with inverse $\tau$. Let $\sigma=\eta^{-1}s$ be a splitting for $\sigma$ and let $(\ep^*)^{-1}t^*$ be a splitting for $\tau^*$. We have $st=\eta\ep$.
\end{cor}
\begin{pf} Since $\sigma\tau=1$, we must have $\eta^{-1}st\ep^{-1}=1$, and hence $st=\eta\ep$, as claimed.
\end{pf}

Using the splitting theorem, we can also start our induction.

\begin{prop} Let $G$ be a UC-proper, supersoluble group. If $\sigma$ is a unit of length $1$, then $\sigma$ is trivial.
\end{prop}
\begin{pf} Since $G$ is UC-proper, let $N$ be a normal subgroup whose quotient is infinite dihedral, generated by $Nx$ and $Ny$. Since $\sigma$ has length $1$, it lives either in $\gen{N,x}$ or $\gen{N,y}$, both of which are subgroups of infinite index in $G$, and hence support no non-trivial units. This proves the result.
\end{pf}

As a corollary, we get an important piece of information.

\begin{cor}\label{etanotone} Let $G$ be a torsion-free supersoluble group, and let $\sigma$ be a unit of $KG$, of length $n$ beginning in $x$. Let
\[ \sigma=\prod_{i=1}^n (\alpha_i+\beta_i\gamma_i)\eta^{-1}\]
be a splitting for $\sigma$. If $\eta$ is a unit then $\sigma$ is a trivial unit.
\end{cor}
\begin{pf} Since $\eta=1$, this implies that $\prod_{i=1}^n (\alpha_i+\beta_i\gamma_i)\tau=1$, where $\tau=\sigma^{-1}$; then $\alpha_n+\beta_n\gamma_n$ is a unit, and since there are no non-trivial length-$1$ units, we have a contradiction.
\end{pf}

In turn, this gives us the result for length $2$.

\begin{cor}\label{cor:nolength2} Let $G$ be a UC-proper, supersoluble group. If $\sigma$ is a unit of length $2$ then $\sigma$ is trivial.
\end{cor}
\begin{pf} Let $\sigma=\eta^{-1}s$ be a splitting for $\sigma$. Expanding out $(\alpha_2+\beta_2x)(\alpha_1+\beta_1y)$ (with $\alpha_i$ and $\beta_i$ left-coprime, which we may assume by pulling out their left-gcds), we get
\[ \alpha_2\alpha_1+\alpha_2\beta_1y+\beta_2\alpha_1^xx+\beta_2\beta_1^xxy,\]
where $\alpha^x=x\alpha x^{-1}$. If $p$ is a prime dividing $\alpha_2\alpha_1$, then it either divides $\alpha_2$ or $\alpha_1$; in the former case, it divides both $\beta_2\alpha_1^x$ and $\beta_2\beta_1^x$, and since $\alpha_1^x$ and $\beta_1^x$ are coprime, we get a contradiction to $\alpha_2$ and $\beta_2$ being coprime. Similarly, we get a contradiction if $p\mid \alpha_1$. Hence, in any splitting of length $2$, the left-gcd of all coefficients of words in $I$ is $1$. Now write $\sigma=\eta^{-1}s$, and note that the left-gcd of the coefficients of the words in $I$ is $1$. In order for $\eta^{-1}s$ to lie in $KG$, we therefore have that $\eta$ is a unit, contradicting Corollary \ref{etanotone}; hence there are no length-$2$ units, as claimed.
\end{pf}

It might be thought that this trend will continue; that is, there can never be a non-trivial $\eta$ dividing all of the coefficients in front of the words in $I$, assuming that the splitting is reduced. This is false, as Example \ref{ex:dividesallcoeffs} demonstrates.

\section{The (Passman) Fours Group}

The `simplest' example of a torsion-free group that is not right-orderable was given by Passman, and is the group
\[ \Gamma=\Gen{x,y}{y^{-1}x^2y=x^{-2},\;x^{-1}y^2x=y^{-2}}.\]
For our work we define $z=xy$, $a=x^2$, $b=y^2$ and $c=z^2$. Then $H=\gen{a,b,c}$ is a normal subgroup of $\Gamma$ isomorphic with $\Z\times \Z \times \Z$, and whose quotient is a Klein four group. Also, $N=\gen{a,b}$ is a normal subgroup of $\Gamma$ isomorphic with $\Z\times\Z$, and whose quotient is infinite dihedral. Let $K$ be a field; then any element $\alpha$ of the group algebra $K\Gamma$ may be written as a sum
\[ \alpha=Ax+By+C+Dz,\]
where $A$, $B$, $C$ and $D$ are elements of $KH$. The group algebra $KH$ may be thought of as a Laurent polynomial ring in three variables, with coefficients in $K$, and we will use this approach. The set $\{1,x,y,z\}$ forms a transversal to $H$ in $\Gamma$, and we will use this as a basis of an embedding of $K\Gamma$ into a matrix ring over $KH$. More precisely, let
\[ x_1=1,\quad x_2=x,\quad x_3=y,\quad x_4=xy.\]
Then there is a $K$-algebra embedding
\[ \eta: K\Gamma\to M_4(KH), \qquad \alpha \mapsto \pi_H(x_i\alpha x_j^{-1}),\]
where $\pi_H$ is the restriction map from $K\Gamma$ to $KH$. If $\alpha$ is written as above, then
\[ \eta(\alpha)=\begin{pmatrix}C&A&B&D\\A^xa&C^x&D^xa&B^x\\B^yb&D^ya^{-1}c^{-1}&C^y&A^ya^{-1}bc^{-1}\\D^zc&B^zb^{-1}&A^zb^{-1}c&C^z\end{pmatrix}.\]
(Here, $A^x$ indicates the conjugate of $A$ by $x$, and so on.) We observe that this representation extends naturally to $\eta:(K\Gamma)(KN)^{-1}\hookrightarrow M_4\((KH)(KN)^{-1}\)$. 

\begin{prop} There are exactly three normal subgroups, $N_1=N$, $N_2=\gen{a,c}$, and $N_3=\gen{b,c}$, such that if $\phi:\Gamma\to D_\infty$ is a surjective homomorphism then $\ker\phi=N_i$ for some $i$. Furthermore, there is an automorphism $\psi$ of $\Gamma$ such that $N_i^\psi=N_{i+1}$ (where the indices are taken modulo $3$).
\end{prop}
\begin{pf} Notice that $(z^2)^x=(xyxy)^x=yxyx$, and 
\begin{align*} (xyxy)(yxyx)&=xyx(y^2)xyx
\\ &=xy(y^{-2})x^2yx
\\ &=xy^{-1}yx^{-2}x
\\ &=1,\end{align*}
so that $x$ conjugates $z^2$ to $z^{-2}$. Similarly, it is easy to see that $y$ also conjugates $z^2$ to $z^{-2}$. Therefore any ordered pair from $\{x,y,z\}$ satisfies the relations of the group, and so there are (outer) automorphisms interchanging $(x,y)$ with $(u,v)$, where $u,v\in\{x,y,z\}$. In particular, all of the $N_i$ are $\Aut(\Gamma)$-conjugate.

Firstly, let $G\cong D_\infty$ be generated by elements $\alpha$ and $\beta$. Since every element of $G$ is either of order $2$ or lies inside the cyclic subgroup of index $2$, it cannot be that both $\alpha$ and $\beta$ have infinite order. Also, if one has infinite order, then their product (either $\alpha\beta$ or $\beta\alpha$) has order $2$ as well. This will be important in what follows.

Let $M$ be a normal subgroup of $\Gamma$ such that $\Gamma/M$ is infinite dihedral. Then $\Gamma/M=\gen{Mx,My}$, and so by the previous paragraph exactly two of $Mx$, $My$, and $Mxy$, must have order $2$ in the quotient. Hence $M$ contains one of the $N_i$, say $N_1$. (Since they are all $\Aut(\Gamma)$-conjugate, we may assume that $N_1\leq M$.) Since any quotient of $D_\infty$ is finite, and we know that $\Gamma/N$ is infinite dihedral, we see that $M=N$, as claimed.
\end{pf}

We can see that $\bigcap N_i=1$, and so for a group element $g\in G$, its images modulo each of the quotients $\Gamma/N_i$ is enough to determine it uniquely. Also, since each of the three normal subgroups $N_i$ are $\Aut(\Gamma)$-conjugate, any result proved using one of the length functions is automatically applicable for the other two length functions got in this way.

There are other length functions on the group, obtained by taking two other generators for $\Gamma$ that satisfy the group relations: for example, consider the pair $(x,xyx)$, which together generate $\Gamma$. Then $\gen{x^2,(xyx)^2}=\gen{x^2,y^{-2}}=N$, but here the elements $x$ and $xyx$ are considered to have length $1$, and the element $y=x(xyx)x$ has length $3$.

\medskip

Since we are interested in finding units, we would like a condition for a group ring element to be a unit.

\begin{thm}[Determinant Condition]\label{detcond} Let $K$ be a field and let $\alpha$ be an element of $K\Gamma$. Then $\eta(\alpha)\in K\setminus\{0\}$ if and only if $\alpha$ is a unit.\end{thm}
\begin{pf}
We will use the fact that $\Gamma$ is supersoluble. Assume that $\alpha\in K\Gamma$ is a unit. Then there exist $\alpha_1,\dots,\alpha_n$, $\beta_1,\dots,\beta_n$, $\nu\in K[a^{\pm 1},b^{\pm 1}]$ such that
\[ \alpha=(\alpha_1+\beta_1\gamma_1)\ldots(\alpha_n+\beta_n\gamma_n)\nu^{-1},\]
for some $\gamma_i\in \{x,y\}$ with $L(\gamma_1\ldots\gamma_n)=n$. It is easy to see that for $\gamma_i=x$, we have
\begin{align*} \det \eta(\alpha_i+x\beta_i)&=(\alpha_i\alpha_i^x-\beta_i\beta_i^xa)(\alpha_i^y\alpha_i^z-\beta_i^y\beta_i^za^{-1});
\intertext{similarly,}
\det \eta(\alpha_i+y\beta_i)&=(\alpha_i\alpha_i^y-\beta_i\beta_i^yb)(\alpha_i^x\alpha_i^z-\beta_i^x\beta_i^zb^{-1}).
\intertext{Finally,}
\det \eta(\nu^{-1})&=(\nu^{-1})(\nu^{-1})^x(\nu^{-1})^y(\nu^{-1})^z.\end{align*}

Since $\det \eta(\alpha)=\prod \det \eta(\alpha_i+\gamma_i\beta_i)\det \eta(\nu^{-1})$, we get that $\det \eta(\alpha)$ is invariant under conjugation by $x$, $y$, and $z$. If $\alpha$ is a unit of $K\Gamma$, then $\det \eta(\alpha)$ is a unit of $KH$, which is of the form $\lambda a^ib^jc^k$, for some $\lambda\in K\setminus\{0\}$. Therefore, we see that $\det \eta(\alpha)=\lambda\in K\setminus\{0\}$.

\medskip

Conversely, if $\alpha\in K\Gamma$ has a determinant in $K\setminus\{0\}$, then $\eta(\alpha)^{-1}\in M_4(KH)$; expressing $\eta(\alpha)^{-1}$ via the matrix of co-factors of $\eta(\alpha)$ of $\eta(\alpha)$ shows directly that $\eta(\alpha)^{-1}$ lies in the image of $\eta$, so that $\alpha^{-1}\in K\Gamma$.
\end{pf}


The next result shows that, in the splitting theorem given in Section \ref{splitsec}, the difference between $\sigma\tau$ and the split form $\sigma\prod (\alpha_i+\beta_i\gamma_i)$ is a central element.

\begin{thm} Let $\sigma$ and $\tau$ be elements of $K\Gamma$, and assume that $\sigma\tau=\eta\in KN\setminus\{0\}$. Then there exist $\alpha_1,\dots,\alpha_s$, $\beta_1,\dots,\beta_s\in KN$, $\gamma_1,\dots,\gamma_s\in\{x,y\}$ such that
\[ \sigma(\alpha_1+\beta_1\gamma_1)\ldots(\alpha_s+\beta_s\gamma_s)=\eta\eta',\]
for some $\eta'\in KN\setminus\{0\}$, central in $K\Gamma$.\end{thm}
\begin{pf} Since $\Gamma$ is supersoluble, $K\Gamma$ has no non-trivial zero divisors. Moreover, Steps 1 to 3 of Theorem \ref{splittingthm} hold, so that (in the notation of that theorem) $\ep_1\delta_1\in KN$ with $L(\sigma\delta_1)<L(\sigma)$. For $\nu\in KN$, let $\prod\nu$ denote the element $\nu\nu^x\nu^y\nu^z$, and let $\prod'\nu$ denote the element $\nu^x\nu^y\nu^z$. Observe that if $\nu$ is non-zero, then $\prod\nu$ is a non-zero element of $KN$ central in $K\Gamma$. With $\nu=\ep_1\delta_1$, we then have
\[ \ep_1\delta_1\prod\nolimits'(\ep_1\delta_1)\ep_1\tau=\ep_1\tau\prod(\ep_1\delta_1).\]
Since $\ep_1$ is non-zero, we conclude that
\[ \delta_1\prod\nolimits'(\ep_1\delta_1)\ep_1\tau=\tau\prod(\ep_1\delta_1),\]
so that
\[ \sigma\delta_1\left[\prod\nolimits'(\ep_1\delta_1)\ep_1\tau\right]=\sigma\tau\prod(\ep_1\delta_1)=\eta\prod(\ep_1\delta_1).\]
The result now follows by induction.\end{pf}

\section{Length-$3$ Units in $K\Gamma$}
\label{sec:nolength3units}

This section is devoted to a proof of the following theorem.

\begin{thm}\label{thm:nolength3units} There are no non-trivial units of length $3$ in $K\Gamma$.\end{thm}

Assume that $\sigma,\tau$ are non-trivial units in $K\Gamma$ such that $\sigma\tau=\tau\sigma=1$ with $L(\sigma)=L(\tau)=3$, which without loss of generality we assume to have longest word $xyx$. Let $I$ denote the subset of $W$ lying in the support of $\sigma$. The splitting of $\sigma$ gives
\[ \lambda(\alpha_3+\beta_3x)(\alpha_2+\beta_2y)(\alpha_1+\beta_1x)\tau=\eta,\]
where $\eta\in KN$ is central in $KG$, and $\lambda\in KN$ is chosen so that $(\alpha_i$ and $\beta_i$ are coprime for $i=1,2,3$. Writing the split part as $s$, we have $\lambda s\tau=\eta$, and so (in the localization of $KG$ at $KN$) $\eta^{-1}\lambda s=\sigma$. We claim that $\lambda$ is a factor of $\eta$: if not, then write $\tilde \lambda=\lambda/(\eta,\lambda)$, and note that $\sigma$ must therefore have the form $\tilde \lambda\sigma'$ for some $\sigma'\in KG$. The left-gcds of the coefficients of the words in $I$ all have $\tilde \lambda$ as a common factor, so by Proposition \ref{prop:leftgcdis1}, $\tilde \lambda$ is a unit. Hence $\lambda\mid \eta$, as claimed.

Write $\tilde\eta=\eta/\lambda$, so that 
\[ \sigma=\tilde\eta^{-1}(\alpha_3+\beta_3x)(\alpha_2+\beta_2y)(\alpha_1+\beta_1x).\]
Define $D_1=\alpha_1\alpha_1^x-\beta_1\beta_1^xa$, $D_2=\alpha_2\alpha_2^y-\beta_2\beta_2^yb$ and $D_3=\alpha_3\alpha_3^x-\beta_3\beta_3^xa$. By direct computation, the element
\[ \(\frac{\alpha_3^x}{D_3}-\frac{\beta_3}{D_3}x\)\(\frac{\alpha_2^y}{D_2}-\frac{\beta_2}{D_2}y\)\(\frac{\alpha_1^x}{D_1}-\frac{\beta_1}{D_1}x\)\tilde\eta\]
is an inverse for $\sigma$ in $(KN)^{-1}(KG)(KN)^{-1}$, and hence by uniqueness of inverses this element is $\tau$.


The following table records the coefficients in front of the words when one expands out the product $s$ of the linear terms in $\sigma$.
\begin{center}\begin{tabular}{|c|l|}
\hline Word & Coefficient
\\ \hline $xyx$ & $\beta_3\beta_2^x\beta_1^{yx}$
\\ $yx$ & $\alpha_3\beta_2\beta_1^y$
\\ $xy$ & $\beta_3\beta_2^x\alpha_1^{yx}$
\\ $y$ & $\alpha_3\beta_2\alpha_1^y$
\\ $x$ & $\alpha_3\alpha_2\beta_1+\beta_3\alpha_2^x\alpha_1^x$
\\ $1$ & $\alpha_3\alpha_2\alpha_1+\beta_3\alpha_2^x\beta_1^xx^2$
\\ \hline
\end{tabular}\end{center}
Since this expanded form is a unit in $K\Gamma$, $\tilde\eta$ must be a factor of each of the coefficients in this table. This allows us to prove the following proposition.

\begin{prop}\label{prop:allprimesdividebeta2} Let $p$ be a prime that divides each of the coefficients of the words in $I$. We have that $p\mid\beta_2,\beta_2^x$, and $p\nmid \alpha_2,\alpha_2^x,\alpha_3,\beta_3$. In particular, $\tilde\eta\mid \beta_2$ and $\tilde\eta\mid\beta_2^x$.
\end{prop}
\begin{pf}We proceed in stages, reducing the problem one step at a time.

\medskip

\noindent\textbf{Step 1}: \textit{Either $p\mid \alpha_3$ or $p\mid\beta_2$, and either $p\mid \beta_3$ or $p\mid\beta_2^x$.} Considering the coefficients of $yx$ and $y$, we see that $p$ divides both $\alpha_3\beta_2\beta_1^y$ and $\alpha_3\beta_2\alpha_1^y$. As $p$ cannot divide both $\alpha_1^y$ and $\beta_1^y$, we must have that either $p\mid\alpha_3$ or $p\mid\beta_2$. Similarly, considering the coefficients of $xyx$ and $xy$, we see that $p$ divides both $\beta_3\beta_2^x\beta_1^{yx}$ and $\beta_3\beta_2^x\alpha_1^{yx}$, so divides either $\beta_3$ or $\beta_2^x$, proving the claim.

Notice that since $p$ cannot divide both $\alpha_3$ and $\beta_3$, if $p\mid\alpha_3$ then $p\mid\beta_2^x$, and similarly if $p\mid\beta_3$ then $p\mid\beta_2$.

\medskip

\noindent\textbf{Step 2}: \textit{$p\nmid \alpha_3$, and so $p\mid\beta_2$.} Suppose that $p\mid\alpha_3$. Since this means that $p\mid\beta_2^x$, we must have that $p\nmid\alpha_2^x$. Considering the coefficients of $x$ and $1$, we see that $p$ divides the first expression in both cases, and so $p\mid\beta_3\alpha_2^x\alpha_1^x,\beta_3\alpha_2^x\beta_1^xa$. This yields a contradiction, since $p\nmid\beta_3$ and $p\nmid\alpha_2^x$. Hence $p\nmid\alpha_3$, so by Step 1, $p\mid\beta_2$.

\medskip

\noindent\textbf{Step 3}: \textit{$p\nmid \beta_3$, and so $p\mid\beta_2^x$.} Suppose that $p\mid\beta_3$. Since this means that $p\mid\beta_2$, we must have that $p\nmid\alpha_2$. Considering the coefficients of $x$ and $1$, we see that $p$ divides the second expression in both cases, and so $p\mid\alpha_3\alpha_2\beta_1,\alpha_3\alpha_2\alpha_1$. This yields a contradiction, since $p\nmid\alpha_3$ and $p\nmid\alpha_2$. Hence $p\nmid\beta_3$, so by Step 1, $p\mid\beta_2^x$. This completes the proof, since $p\nmid\alpha_2,\alpha_2^x$ now.
\end{pf}

\begin{lem}\label{lem:zeroifunit} Let $\alpha$ and $\beta$ be elements of $KN$, and suppose that $\alpha\alpha^y-\beta\beta^yb$ is a unit. Then either $\alpha=0$ or $\beta=0$.
\end{lem}
\begin{pf} By extending $K$ if necessary, we assume that $K$ is infinite. If $u$ is a unit in $KN$, then we may specialize $a$ to be any element of $K$ and the specialization of $u$ remains a unit. Hence specializing $a=k\in K$ yields a polynomial $(\bar\alpha)^2-(\bar\beta)^2b=b^i$. Suppose that both $\bar\alpha$ and $\bar\beta$ are non-zero. Notice that the highest and lowest powers of $b$ in $(\bar\alpha)^2$ are of even degree, and the highest and lowest powers of $b$ in $(\bar\beta)^2b$ are of odd degree. Hence either all of the powers of $b$ in $(\bar\alpha)^2$ are lower than some power of $(\bar\beta)^2$ or vice versa, and similarly either all of the powers of $b$ in $(\bar\alpha)^2$ are larger than some power of $(\bar\beta)^2$ or vice versa. Thus there must be at least two different powers of $b$ present in $(\bar\alpha)^2-(\bar\beta)^2b$, and hence it is not a unit. Thus either $\bar\alpha$ or $\bar\beta$ is zero for the specialization $a=k$. However, if $K$ is infinite then there are infinitely many choices of specialization, but $\bar\alpha$ and $\bar\beta$ can only be zero for finitely many choices of specialization. Thus either $\alpha=0$ or $\beta=0$, as claimed.
\end{pf}

We now embark on the proof of Theorem \ref{thm:nolength3units}, and proceed in stages.

\medskip

\noindent\textbf{Step 1}: \textit{$(D_2,\beta_2)=(\beta_2,\alpha_2^y)$}. Let $A_1=(D_2,\beta_2)$ and $A_2=(\alpha_2^y,\beta_2)$. Since $A_1$ divides both $D_2$ and $\beta_2$, it must divide $\alpha_2\alpha_2^y$; however, since $\beta_2$ and $\alpha_2$ are coprime, $A_1\mid \alpha_2^y$, so that $A_1\mid A_2$. Conversely, $A_2$ divides both $\beta_2$ and $\alpha_2^y$, hence it divides $D_2$; thus $A_2\mid A_1$, so that $A_1=A_2$.

\medskip

The next two stages involve understanding the quotient $D_2'=D_2/(\alpha_2^y,\beta_2)$.

\medskip

\noindent\textbf{Step 2}: \textit{If $p\mid D_2'$ then $p^y\nmid D_2'$}. Write
\[ \tau=\(\frac{\alpha_3^x}{D_3}-\frac{\beta_3}{D_3}x\)\(\frac{\alpha_2^y}{D_2}-\frac{\beta_2}{D_2}y\)\(\frac{\alpha_1^x}{D_1}-\frac{\beta_1}{D_1}x\)\tilde\eta.\]
Let $\alpha'_2=\alpha_2/(\alpha_2,\beta_2^y)$ and $\beta_2'=\beta_2/(\alpha_2^y,\beta_2)$, so that we have
\[ \tau=\(\frac{\alpha_3^x}{D_3}-\frac{\beta_3}{D_3}x\)\(\frac{\alpha_2'^y}{D_2'}-\frac{\beta'_2}{D_2'}y\)\(\frac{\alpha_1^x}{D_1}-\frac{\beta_1}{D_1}x\)\tilde\eta.\]

Applying the regular representation and taking determinants, we get that the expression
\[ \omega=\frac{\alpha_2'^y(\alpha_2'^y)^y-\beta_2'\beta_2'^yb}{D_2'D_2'^y}=\frac{\alpha_2'\alpha_2'^y-\beta_2'\beta_2'^yb}{D_2'D_2'^y}\]
is a factor of $\tilde\eta\tilde\eta^x\tilde\eta^y\tilde\eta^z$. We next notice that $(\alpha_2,\beta_2^y)(\alpha_2'\alpha_2'^y-\beta_2'\beta_2'^yb)=D_2'$, so in fact
\[ \omega=\frac{1}{(\alpha_2,\beta_2^y)D_2'^y};\]
hence $D_2'^y\mid \tilde\eta\tilde\eta^x\tilde\eta^y\tilde\eta^z$. Suppose that $p$ is a prime dividing both $D_2'$ and $D_2'^y$. Therefore $p$ divides $\tilde\eta^\gamma$ for some $\gamma$, and clearly either $p$ or $p^y$ divides either $\tilde\eta$ or $\tilde\eta^x$. Hence, replacing $p$ by $p^y$ if necessary, either $p\mid \tilde\eta$ or $p\mid\tilde\eta^x$. However, by Proposition \ref{prop:allprimesdividebeta2}, all primes dividing $\tilde\eta$ divide both $\beta_2$ and $\beta_2^x$, so $p$ divides both $D_2$ and $\beta_2$. Hence $p$ divides $(D_2,\beta_2)$, so does not divide $D_2'$, a contradiction.

Hence there cannot be a prime dividing both $D_2'$ and $D_2'^y$, as required.

\medskip

\noindent\textbf{Step 3}: \textit{$D_2'=(\alpha_2^y,\beta_2)^y$}. Firstly, $D_2=D_2^y$, so since $D_2'\mid D_2$, we see that $D_2'^y\mid D_2$. By Step 2, $D_2'$ and $D_2'^y$ are coprime, so that, since both $D_2'$ and $D_2'^y$ divide $D_2$, we must have $D_2'D_2'^y\mid D_2$. Finally, by construction of $D_2'$, we must have that $D_2'^y\mid (D_2/D_2')=D_2/(D_2,\beta_2)$, so that $D_2'\mid (\alpha_2^y,\beta_2)^y$. To see the converse, notice that $(\alpha_2^y,\beta_2)^y=(\alpha_2,\beta_2^y)$, which must be prime to $(\alpha_2^y,\beta_2)$. Since $D_2=D_2^y$, $(\alpha_2^y,\beta_2)^y\mid D_2$, and it is prime to $(\alpha_2^y,\beta_2)$, hence divides $D_2'$. Thus we get equality, as claimed.

\medskip

We conclude that $D_2=(\alpha_2^y,\beta_2)(\alpha_2,\beta_2^y)$. In particular, $D_2=A_1A_1^y$, and so
\[ \(\alpha_2/A_1^y\)\(\alpha_2^y/A_1\)-\(\beta_2/A_1\)\(\beta_2^y/A_1^y\)b\]
is a unit, with $\alpha=\alpha_2/A_1^y$ and $\beta=\beta_2/A_1$ elements of $KN$. Hence we have that $\alpha\alpha^y-\beta\beta^yb$ is a unit, so that either $\alpha$ or $\beta$ is zero, by Lemma \ref{lem:zeroifunit}. Clearly $\beta_2\neq 0$, else this element does not have length $3$. However, if $\alpha_2=0$ then $\beta_2$ is a (trivial) unit of $KN$, as $(\alpha_2,\beta_2)=1$. Therefore $\tilde\eta$ is a trivial unit of $KN$, using $\tilde\eta\mid \beta_2$, so that $\sigma =(\alpha_3+\beta_3x)(\alpha_2+\beta_2y)(\alpha_1+\beta_1x)$. Hence each linear factor is a unit in $K\Gamma$, and therefore trivial by the length-one case. This implies that $\sigma$ is a trivial unit of $K\Gamma$, contrary to assumption. This contradiction proves that $\sigma$ is not a unit, and so concludes the proof of Theorem \ref{thm:nolength3units}.

\bigskip

\begin{example}\label{ex:dividesallcoeffs} In Section \ref{sec:conseqsplitting} we proved that for a putative non-trivial unit $\sigma$ of length $2$, the left-gcds of all coefficients in $I$ was $1$, so that $\eta=1$, and $\sigma$ cannot exist (Corollary \ref{cor:nolength2}). A similar strategy will not work for length $3$ units, since it is possible to find $\alpha_i$ and $\beta_i$ for $i=1,2,3$ such that the left-gcd of the coefficients of all words in $I$ is not a unit.

Choose 
\[ \alpha_1=\alpha_2=\alpha_3=\beta_3=1,\quad \beta_1=-a,\quad \beta_2=1-a.\]
We have
\[ (1+x)(1+(1-a)y)(1-ax)=(a-1)\(a^{-1}xyx+a^{-1}yx-xy-y-x+(1+a)\).\]
Of course, this is not a unit, either because of Theorem \ref{thm:nolength3units} or by direct computation.
\end{example}

\section{The Promislow Set}

In \cite{promislow1988}, Promislow constructed a fourteen-element subset $\ms P$ of the Passman fours group $\Gamma$ such that $\ms P\cdot \ms P$ has no unique product. We use the main theorem of the previous section to conclude that it cannot be the support of a unit in $K\Gamma$, for any field $K$.

\begin{thm} Let $K$ be any field, let $\Gamma=\Gen{x,y}{y^{-1}x^2y=x^{-2},\;x^{-1}y^2x=y^{-2}}$ be the Passman fours group, and write $a=x^2$, $b=y^2$, $c=(xy)^2$. Let $\ms P\subset \Gamma$ be the Promislow set
\[\ms P=\ms Ax\cup \ms By\cup \ms C,\]
where
\[ \ms A=\{1,a^{-1},a^{-1}b,b,a^{-1}c^{-1},c\},\quad \ms B=\{1,a,b^{-1},b^{-1}c,c,ab^{-1}c\},\quad \ms C=\{c,c^{-1}\}.\]
There is no unit in $K\Gamma$ whose support is $\ms P$.
\end{thm}
\begin{pf} By Theorem \ref{thm:nolength3units}, $K\Gamma$ has no units of length $3$. Applying the automorphism that fixes $y$ and swaps $x$ and $xy$, (and hence swaps $a$ and $c$, we note that the image of the Promislow set is
\[\ms P'=\ms B'y\cup \ms C'\cup\ms D'xy,\]
where 
\[ \ms B'=\{1,c,b^{-1},b^{-1}a,a,cb^{-1}a\},\quad \ms C'=\{a,a^{-1}\},\quad\ms D'=\{1,c^{-1},c^{-1}b,b,c^{-1}a^{-1},a\}.\]
It is clear all elements of this set not involving $c$ have length at most $2$, since they are of the form $\alpha$, $\alpha y$, and $\alpha xy$ for some $\alpha\in KN$, where $N=\gen{a,b}$. The remaining elements are of the form $\alpha cy$ and $\alpha c^{-1}xy$ for some $\alpha\in KN$. In the former case, this has length $3$ as it is of the form $\alpha'xyx$, and in the latter case it has length $2$, since $c^{-1}xy=y^{-1}x^{-1}=ab^{-1}yx$. Hence any element of $K\Gamma$ with support $\ms P'$ has length $3$, so is not a non-trivial unit of $K\Gamma$, as required.
\end{pf}

\section{The Higher-Length Case}

Let $\sigma$ be a non-trivial unit, and let $\sigma=\eta^{-1}s$ be a splitting for $\sigma$. As we have mentioned, $\eta$ must divide the coefficients of the words in $I$. Proposition \ref{prop:allprimesdividebeta2} proved that, if $L(\sigma)=3$, then all primes dividing $\eta$ divide $\beta_2$ and $\beta_2^x$. When the length of $\sigma$ is greater than $3$, however, there is no unique collection of the $\alpha_i$ and $\beta_i$ that a prime dividing $\eta$ need divide.

Let $n$ be a natural number, and expand the expression
\[ (\alpha_n+\beta_n\gamma_n)(\alpha_{n-1}+\beta_{n-1}\gamma_{n-1})\ldots(\alpha_1+\beta_1\gamma_1),\]
where $\gamma_i\in \{x,y\}$ and $\gamma_i\neq \gamma_{i+1}$. The coefficients in front of the words in $W$ will be denoted by $V_{n,x}$ if $\gamma_n=x$ and $V_{n,y}$ if $\gamma_n=y$.

A collection $M$ of conjugates of the $\alpha_i$ and $\beta_i$ is called a \emph{consistent chain} for $V_{n,x}$ (and similarly for $V_{n,y}$) if
\begin{enumerate}
\item whenever $v$ is an element of $V_{n,x}$ and all but one of the terms in $v$ contain an element of $M$, then all terms in $v$ contain an element of $M$, and
\item whenever $\alpha_i^\gamma$ lies in $M$, $\beta_i^\gamma$ does not, and whenever $\beta_i^\gamma$ lies in $M$, $\alpha_i^\gamma$ does not.
\end{enumerate}
A consistent chain is a set $R$ such that if $p$ is a prime dividing all elements of $V_{n,x}$, then $p$ can divide all elements of $R$ without dividing all but one of the terms in any element of $V_{n,x}$; if $p$ divided all but one of the terms in an element of $V_{n,x}$, then $p$ must divide the last, and so divides one of the $\beta_i^\times$ or $\alpha_i^\times$ (where $\times$ is one of $x$, $y$, $z$ or nothing). We illustrate the concept of a consistent chain with an example.

\begin{example} In Section \ref{sec:nolength3units} we described in a table the set $V_{3,x}$. A consistent chain for these is, for example, the set $\{\beta_2,\beta_2^x\}$, or $\{\beta_2,\beta_2^x,\alpha_1,\beta_1^x\}$. Proposition \ref{prop:allprimesdividebeta2} proves that all consistent chains contain $\{\beta_2,\beta_2^x\}$ as a subset, and no consistent chain contains either $\alpha_3$ or $\beta_3$.
\end{example}

In this section we give a recursive description of the `minimal' consistent chains for $V_{n,x}$ and $V_{n,y}$, minimal in the sense that any consistent chain for $V_{n,x}$ contains a minimal one as a subset. Define $U_{n,x}$ to contain the elements $\beta_{n-1}$, $\beta_{n-2}^y$, $\beta_{n-3}^z$, $\beta_{n-4}^x$, $\beta_{n-5}$, and repeating this sequence until the appropriate conjugate of $\beta_2$, and $U_{n,y}$ to be the same sequence with $y$ swapped with $x$.

In the proof of this theorem we will need to understand certain elements of $V_{n,x}$, and so it will help to have the following small-length examples as a guide.

\begin{center}\begin{tabular}{|c|cc|c|cc|}
\hline Length & Word & Coefficient & Length & Word & Coefficient
\\ \hline  \multirow{4}{*}{$4$}&$xyxy$ & $\beta_4\beta_3^x\beta_2^{yx}\beta_1^{xyx}$ & \multirow{4}{*}{$5$}&$xyxyx$ & $\beta_5\beta_4^x\beta_3^{yx}\beta_2^{xyx}\beta_1^{yxyx}$
\\ &$yxy$ & $\alpha_4\beta_3\beta_2^y\beta_1^{xy}$ & &$yxyx$ & $\alpha_5\beta_4\beta_3^y\beta_2^{xy}\beta_1^{yxy}$
\\ &$xyx$ & $\beta_4\beta_3^x\beta_2^{yx}\alpha_1^{xyx}$ & &$xyxy$ & $\beta_5\beta_4^x\beta_3^{yx}\beta_2^{xyx}\alpha_1^{xyxy}$
\\ &$yx$ & $\alpha_4\beta_3\beta_2^y\alpha_1^{xy}$ & &$yxy$ & $\alpha_5\beta_4\beta_3^y\beta_2^{xy}\alpha_1^{yxy}$
\\ \hline
\end{tabular}\end{center}

\begin{thm} Let $n\geq 3$ be an integer. The minimal consistent chains $M_{n,x}$ for $V_{n,x}$ are all pairs $\{\lambda,\mu\}$, with $\lambda$ and $\mu^x$ appearing in the list $U_{n,x}$, together with the minimal consistent chains for $V_{n-1,y}$ (i.e., $\{R\cup\{\beta_n\}:R\in M_{n-1,y}\}$) together with, and those for $V_{n-1,y}$ conjugated by $x$ together with $\alpha_n$ (i.e., $\{R^x\cup\{\alpha_n\}:R\in M_{n-1,y}\}$). The minimal consistent chains $M_{n,y}$ for $V_{n,y}$ are the same, with $x$ and $y$ swapped.
\end{thm}
\begin{pf}
Without loss of generality, assume that $\gamma_n=x$. Let $R$ denote a consistent chain, and suppose firstly that $\beta_n\in R$. We may remove all of the terms from $V_{n,x}$ that start with $\beta_n$ to get a set $V_{n,x}^*$, and by considering
\begin{equation}\label{eq:factoredform}(\alpha_n+\beta_n\gamma_n)(\alpha_{n-1}+\beta_{n-1}\gamma_{n-1})\ldots(\alpha_1+\beta_1\gamma_1),\end{equation}
we clearly see that
\[ V_{n,x}^*=\{\alpha_nw:w\in V_{n-1,y}\}.\]
Since $\alpha_n\notin R$, we may remove the $\alpha_n$ from the start of the words in $V_{n,x}^*$, and so $R\setminus\{\beta_n\}$ must be a consistent chain for $V_{n-1,y}$, as $R$ is a consistent chain for $V_{n,x}$. This case is covered in the theorem, so we may assume that $\beta_n$ does not lie in $R$.

Similarly, suppose that $\alpha_n$ lies in $R$. In this case we may remove all of the terms from $V_{n,x}$ that start with $\alpha_n$ to get a set $V_{n,x}^*$, and we see that 
\[ V_{n,x}^*=\{\beta_nw^x:w\in V_{n-1,y}\}.\]
As above, the elements $R\setminus\{\alpha_n\}$ conjugated by $x$ form a consistent chain for $V_{n-1,y}$, and this case is also covered in the theorem. Hence we may assume that neither $\alpha_n$ nor $\beta_n$ lie in $R$.

We now note that, when expanding (\ref{eq:factoredform}), there are four elements of $V_{n,x}$ that are monomials, namely the coefficients of the words of lengths $n$, $n-1$, and the word of length $n-2$ starting in $y$: two of these words start with $x$, and two start with $y$. If $a_1$ and $a_2$ are the two monomial coefficients of the words starting in $x$, then 
\[ a_1=\beta_n\beta_{n-1}^x\beta_{n-2}^z\ldots\beta_1^\times,\quad a_2=\beta_n\beta_{n-1}^x\beta_{n-2}^z\ldots\alpha_1^\times\]
(where $\times$ is one of $x$, $y$, $z$, or nothing, and for the rest of the proof will also denote one of these four). Since $a_1$ and $a_2$ differ only in the last element, if $R$ is a consistent chain then $R$ must contain at least one of the terms $\beta_i^\times$ for $1<i<n$. Similarly, if $b_1$ and $b_2$ denote the two monomial coefficients of the words starting in $y$, then
\[ b_1=\alpha_n\beta_{n-1}\beta_{n-2}^y\ldots\beta_1^\times,\quad b_2=\alpha_n\beta_{n-1}\beta_{n-2}^y\ldots\alpha_1^\times.\]
Again, $b_1$ and $b_2$ differ only in the last element, so if $R$ is a consistent chain then $R$ must contain at least one of the terms $\beta_i^\times$ for $1<i<n$. It remains to note that the middle $\beta_i^\times$ of the $b_i$ are $U_{n,x}$, and the middle $\beta_i^\times$ of the $a_i$ are the elements of $U_{n,x}$ conjugated by $x$. Thus $R$ contains $\{\lambda,\mu\}$, where $\lambda,\mu^x\in U_{n,x}$, as claimed by the theorem.
\end{pf}

If $\sigma$ is a non-trivial unit of length $n$, starting in $x$, then the $\eta$ obtained from the split form is non-trivial, and any prime $p$ dividing $\eta$ must divide each of the elements of $V_{n,x}$. Hence $p$ must be a factor of every element of a minimal consistent chain $R$.

If $n=3$ then there is only one minimal consistent chain for $V_{3,x}$, namely $\{\beta_2,\beta_2^x\}$. For $n=4$ there are more minimal consistent chains for $V_{4,x}$, namely
\[ \{\beta_3,\beta_3^x\},\;\{\beta_3,\beta_2^z\},\;\{\beta_2^y,\beta_3^x\},\;\{\beta_2^y,\beta_2^z\},\;\{\beta_4,\beta_2,\beta_2^y\},\;\{\alpha_4,\beta_2^x,\beta_2^z\},\]
and the minimal consistent chains for $V_{4,y}$ are
\[ \{\beta_3,\beta_3^y\},\;\{\beta_3,\beta_2^z\},\;\{\beta_2^x,\beta_3^y\},\;\{\beta_2^x,\beta_2^z\},\;\{\beta_4,\beta_2,\beta_2^x\},\;\{\alpha_4,\beta_2^y,\beta_2^z\}.\]
However, some of these are related by applying automorphisms. Denote by $c_x$ conjugation by $x$, $c_y$ conjugation by $y$, $\phi$ the automorphism interchanging $x$ and $y$, and $*$ for the usual anti-automorphism. Applying the anti-automorphism $*$ sends units to units, and applies the map
\[ (\alpha_4+\beta_4x)(\alpha_3+\beta_3y)(\alpha_2+\beta_2x)(\alpha_1+\beta_1y)\mapsto(\alpha_1^z+\beta_1^xy)(\alpha_2^z+\beta_2^yx)(\alpha_3^z+\beta_3^xy)(\alpha_4^z+\beta_4^yx).\]
Applying these automorphisms of $\Gamma$ permutes the minimal consistent chains. For example, suppose that $p$ divides $\{\beta_3^x,\beta_2^y\}$: conjugating by $x$ yields a prime $p$ dividing $\sigma^x$ that divides $\{\beta_3,\beta_2^z\}$. In fact, using these automorphisms we can divide the minimal consistent chains into two collections.
\[ \begin{diagram}\node{\{\beta_3,\beta_3^x\}}\arrow{s,l}{\phi}\arrow{e,t}{\ast}\node{\{\beta_2,\beta_2^y\}}\arrow{s,r}{\phi}\arrow{e,t}{c_x}\node{\{\beta_2^x,\beta_2^z\}}
\\ \node{\{\beta_3,\beta_3^y\}}\node{\{\beta_2,\beta_2^x\}}\arrow{e,t}{c_x}\node{\{\beta_y,\beta_2^z\}}
\\ \node{\{\beta^x_3,\beta_2^y\}}\arrow{e,b}{c_x}\node{\{\beta_3,\beta_2^z\}}\arrow{e,b}{c_y}\node{\{\beta_2^x,\beta_3^y\}}
\end{diagram}\]
(Note that not all arrows are on this diagram.) Suppose that one can prove that there is no unit $\sigma$ of length $4$ and prime $p\mid \eta$ such that $p$ divides $\{\beta_3,\beta_3^x\}$ or $\{\beta_3,\beta_2^z\}$. By the diagram above, applying automorphisms of $\Gamma$ proves that there are no consistent chains that $\eta$ can divide, so $\eta$ is trivial. This allows us to drastically reduce the number of minimal consistent chains that need to be considered when proving that no non-trivial units of length $n$ exist.

\bibliography{references}

\end{document}